\documentclass[12pt]{article}
\usepackage[left=1in,right=1in,top=1in,bottom=1in]{geometry}
\usepackage{amsmath}
\usepackage{amssymb}
\usepackage{graphicx}
\usepackage{listings}
\usepackage[shortlabels]{enumitem}
\usepackage{hyperref}
\usepackage{algorithmicx,algorithm}
\usepackage{algpseudocode}
\usepackage{float}
\usepackage[super]{nth}
\usepackage{authblk}
\usepackage{caption}
\usepackage{subcaption}

\usepackage{etoolbox}
\usepackage{url}
\usepackage[most]{tcolorbox}
\usepackage{xcolor}
\usepackage{bm}
\usepackage[thinc]{esdiff}
\usepackage{pdfpages}
\usepackage{mathtools}

\DeclarePairedDelimiter\abs{\lvert}{\rvert}%
\DeclarePairedDelimiter\norm{\lVert}{\rVert}%

\DeclareMathOperator{\tr}{tr}

\makeatletter
\let\oldabs\abs
\def\abs{\@ifstar{\oldabs}{\oldabs*}}
\let\oldnorm\norm
\def\norm{\@ifstar{\oldnorm}{\oldnorm*}}
\makeatother

\graphicspath{{images/}}

\newtheorem{theorem}                {Theorem}
\newtheorem{lemma}      [theorem]   {Lemma}

\newtheorem{prop}       [theorem]   {Proposition}

\def\endproof{\hfill $\rule{1.3ex}{1.3ex}$ \endtrivlist}

\providetoggle{solution}
\settoggle{solution}{false}

\def\N{\mathcal{N}}

\def\C{\mathcal{C}}

\def\R{\mathbb{R}}
\def\E{\mathbb{E}}

\def\W{\mathcal{W}}
\def\O{\mathcal{O}}

\usepackage{setspace}
\begin{document}

\title{Lower Bounds and a Near-Optimal Shrinkage Estimator for Least Squares using Random Projections}
\author{Srivatsan Sridhar} \author{Mert Pilanci} \author{Ayfer \"Ozg\"ur} \affil{Stanford University \\ \{svatsan, pilanci, aozgur\}@stanford.edu}

\date{}
\maketitle

\begin{abstract} 
     In this work, we consider the deterministic optimization using random projections as a statistical estimation problem, where the squared distance between the predictions from the estimator and the true solution is the error metric. In approximately solving a large scale least squares problem using Gaussian sketches, we show that the sketched solution has a conditional Gaussian distribution with the true solution as its mean. Firstly, tight worst case error lower bounds with explicit constants are derived for any estimator using the Gaussian sketch, and the classical sketching is shown to be the optimal unbiased estimator. For biased estimators, the lower bound also incorporates prior knowledge about the true solution.
    Secondly, we use the James-Stein estimator to derive an improved estimator for the least squares solution using the Gaussian sketch. An upper bound on the expected error of this estimator is derived, which is smaller than the error of the classical Gaussian sketch solution for any given data. The upper and lower bounds match when the SNR of the true solution is known to be small and the data matrix is well conditioned. Empirically, this estimator achieves smaller error on simulated and real datasets, and works for other common sketching methods as well.
\end{abstract}

\section{Introduction}
One of the simplest and most ubiquitous of optimization problems is the least squares problem. The least squares problem applies to several common models such as linear regression, maximum likelihood estimation, and loss functions in neural networks. The unconstrained least squares problem can be defined as
\begin{equation*} x_{LS} = \text{arg}\min\limits_{x\in\R^d} \, \norm{Ax - y}_2^2
\end{equation*}
The data matrix $A\in\R^{n\times d}$ and the data vector $y\in\R^n$ are fixed and known, and $x_{LS}$ is the least squares solution. In this work, we consider $n>d$. This problem has a unique analytical solution $x_{LS}=A^\dagger y$ where $A^\dagger=(A^TA)^{-1}A^T$ is the pseudoinverse of $A$, whenever $rank(A)=d$. This analytical solution can be computed efficiently using singular value decomposition, or Cholesky or QR decomposition, all of which require $\O(nd^2)$ time.

Often in modern datasets, we deal with very large scale data, i.e. $n$ and $d$ are very large, making the computation of the true least squares solution very expensive. We can approximate the least squares solution by projecting $A$ and $y$ onto a lower dimension $m<n$, and then solving the least squares problem using the projected data. This method is known as sketching, and can be defined as
\begin{equation*}
    \widehat{x} = \text{arg}\min\limits_{x\in\R^d} \, \norm{SAx - Sy}_2^2
\end{equation*}
where $S\in\R^{m\times n}$ is a random projection matrix. $m$ is called the sketch size. This sketched problem can be solved in $\O(md^2)$ time, plus the time to compute the sketch $SA$. Apart from reducing computation time, sketching can also be useful to provide differential privacy \cite{blocki,showkatbakhsh,zhou2009}. A reasonable sketch size must be larger than $d$, and grows proportionally to $d$.

The error of the sketched estimator can be measured as the distance between the predictions due to the estimator and the true solution $\norm{A(\widehat x-x_{LS})}_2^2$ (prediction error), which we will consider in this work. It is known that this approximately grows as $d/m$ for several sketching methods \cite{dobriban2019asymptotics, pilanci2014iterative, raskutti2016statistical}.

There are several common sketching methods, discussed in \cite{mahoney2011randomized, sarlos2006, vempala2005random, Woodruff14}. One simple method is forming $SA$ and $Sy$ by sampling $m$ rows of $A$ and $y$ \cite{mahoney2011randomized, sarlos2006, drineas2006, drineas2011faster}.
One can instead use a random matrix $S$ from a data-oblivious distribution. $S_{ij}\overset{\text{i.i.d.}}{\sim}\frac{1}{\sqrt{m}}\N(0,1)$, called the Gaussian sketch, is the canonical example \cite{sarlos2006}. Many other sketching methods, such as the Rademacher sketch, and the Fast Johnson-Lindenstrauss (JL) sketch, tend to behave very similar to the Gaussian sketch.

The Gaussian sketch has many desirable properties. Firstly, the expected error can be exactly computed. Secondly, fundamental to our results, the sketched solution $\widehat x$ has a Gaussian distribution with mean $x_{LS}$, when conditioned on $SA$. The latter property allows us to view the deterministic optimization problem as a task of statistical estimation - estimating the mean $x_{LS}$ from $\widehat x$. To our knowledge, previous studies have not considered the estimation aspect of sketching.

Despite its simplicity and effectiveness, there are very few results about error lower bounds for the Gaussian sketch, and these lower bounds involve constants that are hard to compute \cite{pilanci2014iterative, raskutti2016statistical, pilanci2015randomized}. Thus there are no results about the optimality of the sketching methods. In this paper, we provide two error lower bounds for the Gaussian sketch. The first one which holds for unbiased estimators only, shows that the classical sketching method is the optimal unbiased estimator. The second lower bound holds for general estimators in the presence of prior knowledge about the solution, and is in terms of explicit constants. Further, much of the previous work assumes a statistical observation model on the data, such as $y=Ax^*+e$, where $e$ is Gaussian noise, which is restrictive and may not be valid in many cases. Instead, our lower bounds assume no such model. Interestingly, the statistical nature of the problem arises because of the randomness in the sketching only, and not from the data model.

We also note that standard sketching methods are based on a sketch-and-solve approach (i.e. solving the original problem with the sketched data), which lead to unbiased estimators in the Gaussian sketch. However, this may not be the best approach, e.g., iterative sketching methods perform better \cite{pilanci2014iterative}). For estimating the unknown mean of a high-dimensional Gaussian vector, it is known that a ``shrinkage estimator" (e.g. the James-Stein estimator) has lower variance than the maximum likelihood (ML) estimator. We use this fact to derive an improved estimator for $x_{LS}$, based on the Gaussian sketch. This estimator is biased, however, has lower error than the classical Gaussian sketch, more so when the data is noisy. Furthermore, this estimator achieves our lower bound when the SNR of the true solution is known to be small and the data matrix is well conditioned. Thus, our analysis introduces a bias-variance trade-off in sketching methods to improve the classical sketch-and-solve paradigm.

\subsection{Our Contributions}
We first show two worst case lower bounds on the prediction error of any estimator that can be obtained from the Gaussian sketch. The first one holds for unbiased estimators of $x_{LS}$, and shows that the classical sketching method is the minimum variance unbiased estimator (MVUE) of $x_{LS}$. Although MVUE is a well known concept for statistical estimation under additive Gaussian noise, this is a novel result for sketching the deterministic least squares optimization problem where the only source of randomness is the sketching matrix. The second lower bound holds for biased estimators as well. The second lower bound deals with constrained least squares problems in general, and this allows us to incorporate the effect of prior knowledge about the least squares solution (such as constraints on the norm or SNR). Both the lower bounds are easy to evaluate in terms of properties of the data (and the constraints). To derive the lower bounds, we use the Fisher information, along with the Cram\'er-Rao inequality for the unbiased case, and the van Trees inequality \cite{gill1995} for the general case.

Our second result is an improved estimator for the Gaussian sketch, based on the James-Stein estimator. We provide an approximate form of the estimator which can be evaluated from the Gaussian sketch solution without much additional computation. We derive an upper bound on the prediction error of this estimator, and show that the error is lesser than that of the classical Gaussian sketch for any given data, even more so when the signal-to-noise ratio (SNR) of the data is low. Finally we show that the lower and upper bounds match for large $m$ upto constants in the additive term, when the SNR of the data is known to be small and $A$ is well-conditioned.

Empirically, we show that our estimator achieves lower prediction error for synthetic as well as real data, especially when the SNR is small. We also show that our estimator works equally well for other sketching methods such as sampling and Fast JL sketch. We also show extensions of our results to matrix regression problems.

\subsection{Previous Works}
There are many different sketching methods in literature \cite{mahoney2011randomized, sarlos2006, vempala2005random, Woodruff14}. Among sampling strategies, a simple uniform sampling can perform poorly in the worst case, but sampling based on leverage scores is known to be better \cite{drineas2006}. Other than Gaussian, common sketching matrices are the Rademacher sketch ($S_{ij}\sim\frac{1}{\sqrt{m}}\text{Unif}\{+1,-1\}$), and the Fast Johnson-Lindenstrauss (JL) sketch \cite{drineas2011faster, mahoney2011randomized, boutsidis2009random}.
Multiplying $A$ with a dense matrix $S$ can require $\O(mnd)$ time. But the Fast JL sketch introduced in \cite{ailon2006} (based on a randomized Hadamard transform) can be computed efficiently in $\O(nd\log n)$ time. The Clarkson-Woodruff (CW) sketch \cite{clarkson2017low} can be computed in just $\O(nd)$ time and has performance similar to the Fast JL sketch \cite{Woodruff14, ahfock2017statistical}.

The statistical perspective of sketching is not new, and has been considered in \cite{ahfock2017statistical} for Gaussian, Fast JL and CW sketches, in \cite{li2006very} for sparse Rademacher sketches, and in \cite{ma2015statistical} for leverage score sampling. \cite{dobriban2019asymptotics} and \cite{raskutti2016statistical} analyze asymptotic properties of various sampling and sketching methods. While for Gaussian sketches, it is possible to analyze statistical properties for finite $n,d$, the other sketches are analyzed only in asymptotic regimes. In fact, \cite{ahfock2017statistical} and \cite{li2006very} show that the solutions from Fast JL, CW and sparse Rademacher sketches asymptotically converge to the Gaussian sketch estimator in distribution. Although our analysis of upper and lower bounds is based on the Gaussian sketch, using this result we can argue that our results apply to a broader class of sketching matrices in the asymptotic regime ($n\to\infty$).

Error lower bounds on sketching for convex constrained least squares are shown in \cite{pilanci2014iterative, raskutti2016statistical, pilanci2015randomized}. \cite{pilanci2015randomized} gives lower bounds on the sketch size $m$ for the residual error to be within $(1+\delta)$ times the residual error of the true solution. \cite{pilanci2014iterative} provides a worst case prediction error lower bound using tools from hypothesis testing and Fano's inequality. \cite{raskutti2016statistical} verifies that this lower bound is tight up to constants for Gaussian, Fast JL sketches and certain sampling methods. However that lower bound involves properties of the convex constraint set that are hard to evaluate, and constants that are not tight. In this work, we provide a lower bound for Gaussian sketches, which can be exactly evaluated. This is done by using a different technique for minimax lower bounds using the Fisher information, which was shown in \cite{barnes2019} to give lower bounds with simpler analysis and precise constants. Another difference is that while the lower bound in \cite{pilanci2014iterative} assumes a statistical model on the data, our lower bound is a worst case bound for any given data.

\subsection{Organization}
Section \ref{sec:background} provides the background about properties of the Gaussian sketch and the James-Stein estimator that will be useful to understand our results. Section \ref{sec:lower_bounds} gives the two lower bound results. Section \ref{sec:estimator} gives the improved estimator, and an upper bound on its error, and analyzes the tightness of the upper and lower bounds. Section \ref{sec:empirical} contains empirical results on simulated and real data, on other sketching methods, and on alternate forms of the proposed estimator. Section \ref{sec:conclusion} gives the conclusion and interesting directions opened by this work. Proofs of all theoretical results are given in Section \ref{sec:proofs}.

\section{Background}
\label{sec:background}
\subsection{Gaussian Sketches for Least Squares Optimization}
Consider an unconstrained least squares optimization problem 
\begin{equation} x_{LS} = \text{arg}\min\limits_{x\in\R^d} \, \norm{Ax - y}_2^2
\label{eq:ls}
\end{equation}
where the data matrix $A\in\R^{n\times d}$ and the data vector $y\in\R^n$ are fixed and known. It is assumed that $n>d$. Throughout this project, we will assume that rank$(A)=d$, so that $A^TA$ is invertible. The classical Gaussian sketch is defined as
\begin{equation}
    \widehat{x} = \text{arg}\min\limits_{x\in\R^d} \, \norm{SAx - Sy}_2^2
    \label{eq:sketchls}
\end{equation}
where $A$ and $y$ are replaced by low dimensional projections $SA$ and $Sy$, where $S\in\R^{m\times n}$ has i.i.d. entries $S_{ij}\sim\frac{1}{\sqrt{m}}\N(0,1)$. Note that the sketching matrix $S$ satisfies $\E[S^TS]=I_n$. 
\begin{lemma}
\label{lem:ls}
The classical Gaussian sketch $\widehat{x}$ satisfies
\begin{align}
\label{eq:sketch_solution}
    \widehat{x} &= x_{LS} + (SA)^\dagger Sy^\perp
\end{align}
where $(SA)^\dagger = (A^TS^TSA)^{-1}A^TS^T$ is the pseudoinverse of $SA$ and $y^\perp=y-Ax_{LS}$ such that $A^Ty^\perp=0$. $\norm{y^\perp}_2^2$ is called the residual error. Further, $Sy^\perp$, and the columns of $SA$ are independent multivariate Gaussians, and
\begin{equation}
\label{eq:sketch_dist}
    \widehat{x} \sim \N\left(x_{LS}, \frac{1}{m}\norm{y^\perp}_2^2 (A^TS^TSA)^{-1}\right) \mid SA
\end{equation}
\end{lemma}
Here the notation $\mid SA$ means conditioning on the random variable $SA$. Also Lemma \ref{lem:ls} assumes that $A^TS^TSA$ is invertible, which is true with probability 1, if $A^TA$ is invertible \cite{gupta2018matrix}.


This property allows us to look at the least squares optimization problem as a statistical estimation problem of estimating the mean of a multivariate Gaussian from its sample, conditioned on the matrix $SA$. Note that this property holds irrespective of the distribution of the data, and is a property of the Gaussian sketch method. In particular this means that $\widehat x$ is an unbiased estimator of $x_{LS}$ since $\E[\widehat x]=x_{LS}$.

The expected error of the classical Gaussian sketch can be evaluated exactly. This is known from earlier works as well \cite{ahfock2017statistical}, and is stated below again.

\begin{lemma}
\label{lem:gaussian_sketch}
The prediction error of the classical Gaussian sketch $\widehat x$ is:
\begin{equation}
\label{eq:sketch_error}
    \E\left[\norm{A(\widehat{x}-x_{LS})}_2^2 \right] = \frac{d}{m-d-1} \norm{y^\perp}_2^2
\end{equation}
\end{lemma}

\subsection{James-Stein Shrinkage Estimator}
Suppose we are given $X \sim \mathcal{N}(\theta,I_d)$ where $\theta\in\mathbb{R}^d$ is an unknown mean vector. We want to design an estimator $\widehat{\theta}(X)$ such that it minimizes the expected squared error $R(\widehat{\theta},\theta)=\mathbb{E}\lVert\widehat{\theta}-\theta\rVert_2^2$. The maximum likelihood estimator for $\theta$ is $\widehat{\theta}(X)=X$ which is an unbiased estimator, i.e. $\E\,\widehat\theta=\theta$. However, from \cite{stein1956} and \cite{james1961}, we have that the estimator given by 
\begin{equation}
\label{eq:shrink_I}
    \widehat{\theta}(X) = \left(1-\frac{d-2}{\lVert X\rVert_2^2} \right) X
\end{equation}
results in a lower variance than the maximum likelihood estimator for $d>2$ for any value of $\theta$. In statistical language, $\widehat{\theta}(X)=X$ is an inadmissible estimator and (\ref{eq:shrink_I}) strictly dominates the former. This is a "shrinkage estimator" (as it shrinks the estimate towards 0), now known as the James-Stein estimator.
This estimator has been very useful in statistics, see \cite{efron1977} for example. A simple proof and explanation of this result can be found in \cite{yung1999explaining} and \cite{samworth2012stein}. There are also several improvements and generalizations of this estimator \cite{efron1976, hansen2008generalized, hoff2012shrinkage, shao1994improving} to other shrinkage estimators. Particularly, we look at the generalization to the case where $X\sim\N(\theta,\Sigma)$.
\begin{lemma}
\label{lem:stein}
Suppose $X\sim\N(\theta,\Sigma)$ where $\theta\in\R^d$ is the unknown mean with $d>2$. $\Sigma$ and a single sample $X$ are known. Then consider the estimator for $\theta$ as
\begin{equation}
\label{eq:shrink_Sigma}
    \widehat\theta(X) = \left(1- \frac{d-2}{X^T\Sigma^{-1}X} \right)X
\end{equation}
Then, its expected normalized squared error is given by
\begin{equation}
\label{eq:shrink_error}
    \E[(\widehat\theta - \theta)^T\Sigma^{-1} (\widehat\theta - \theta)] = d - (d-2)^2 \E\left[ \frac{1}{X^T\Sigma^{-1}X} \right]
\end{equation}
\end{lemma}


If we used the estimator $\widehat\theta(X)=X$, then we would have $\E[(X-\theta)^T \Sigma^{-1}(X-\theta)]=d$ and the estimator shown above achieves strictly lesser expected error for any value of $\theta$ when $d>2$. However, this is a biased estimator and the expected squared error will be lower when the true value of $\theta$ is close to zero, as compared to $\theta$ far from zero. More generally, we can also shrink the estimator towards any other point $\theta_0$ and derive an estimator $\theta_0 + \left(1- (d-2)/X^T\Sigma^{-1}X \right)(X-\theta_0)$, so that the error is smaller when $\theta$ is close to $\theta_0$.

The above section motivates us to use the shrinkage estimator to estimate the true least squares solution $x_{LS}$ from the solution $\widehat x$ of one sketched least squares problem.

\section{Error Lower Bounds for Gaussian Sketch}
\label{sec:lower_bounds}
This section presents lower bounds for estimators derived from the Gaussian sketch. The following two theorems give prediction error lower bounds for any estimator of $x_{LS}$ that uses only the Gaussian sketched data $SA$ and $Sy$.  Theorem \ref{theo:unbiased_lower} is a lower bound for unbiased estimators only while Theorem \ref{theo:lower} holds for biased estimators as well. All the proofs are given in Section \ref{sec:proofs}.

\subsection{Lower Bound on Error for Unbiased Estimators}

\begin{theorem}
\label{theo:unbiased_lower}
For any unbiased estimator $\widetilde x$ of $x_{LS}$, (i.e. $\E[\widetilde x]=x_{LS}$) obtained from the Gaussian sketched data $SA$ and $Sy$,
\begin{equation}
\label{eq:unbiased_lower}
    \E\left[\norm{A(\widetilde{x}-x_{LS})}_2^2 \right] \geq \frac{d}{m-d-1}\norm{y^\perp}_2^2
\end{equation}
\end{theorem}
From Lemma \ref{lem:gaussian_sketch}, the classical Gaussian sketch achieves the lower bound for unbiased estimators and is the minimum variance unbiased estimator of $x_{LS}$ using only the Gaussian sketched data $SA$ and $Sy$. Thus we have the tight lower bound
\begin{equation}
\inf_{\widetilde x(SA,Sy) : \E\left[\widetilde x\right]=x_{LS}} \E\left[\norm{A(\widetilde{x}-x_{LS})}_2^2 \right] = \frac{d}{m-d-1}\norm{y^\perp}_2^2\,,
\end{equation}
where the infimum ranges over all unbiased estimators of $x_{LS}$ that are functions of $(SA,Sy)$.
However, it can be seen ahead that when biased estimators are allowed, the classical Gaussian sketch is sub-optimal.

\subsection{A More General Lower Bound}
For the purpose of this theorem, consider a constrained least squares optimization problem
\begin{equation}
    x_{LS} = \text{arg}\min\limits_{x\in\mathcal{C}} \, \norm{Ax - y}_2^2
    \label{eq:ls_cons}
\end{equation}
where $\mathcal{C}$ is a convex constraint set. For the unconstrained problem, set $\mathcal{C}=\R^d$.
\begin{theorem}
\label{theo:lower}
Suppose there exists $B>0$ such that 
$[-B,B]^d\subseteq\mathcal{C}$. Then for any estimator $\widetilde{x}$ for $x_{LS}$ obtained from the Gaussian sketched data $SA$ and $Sy$,
\begin{equation}
\label{eq:lower_B}
    \E\left[\norm{A(\widetilde{x}-x_{LS})}_2^2 \right] \geq \frac{d}{m} \norm{y^\perp}_2^2 \left(1- \frac{\pi^2\norm{y^\perp}_2^2} {mB^2\sigma_{min}(A^TA)} \right)
\end{equation}
where $\sigma_{min}(A^TA)$ is the minimum eigenvalue of $A^TA$. In particular, for unconstrained least squares optimization,
\begin{equation}
\label{eq:lower_infinity}
   \E\left[\norm{A(\widetilde{x}-x_{LS})}_2^2 \right] \geq  \frac{d}{m} \norm{y^\perp}_2^2
\end{equation}
\end{theorem}
The proof will use a method based on the van Trees inequality which uses the Fisher information to lower bound the error \cite{gill1995}. This method has been inspired from \cite{barnes2019} and it provides a neat lower bound in the constrained optimization case too. The lower bound in (\ref{eq:lower_infinity}) must be seen as a worst case lower bound (since it holds for all possible values of $x_{LS}$). In practice, if we wish to impose any regularity conditions on the solution (such as small norm) based on prior knowledge, then we would have a constraint with $B<\infty$ and thus (\ref{eq:lower_B}) would hold which can give us a lower error lower bound.

\section{James-Stein Estimator for Least Squares}
\label{sec:estimator}
\subsection{The Improved Estimator}
Using the Gaussian distribution in (\ref{eq:sketch_dist}) with $\theta=x_{LS}$ and $\Sigma=\frac{1}{m}\norm{y^\perp}_2^2 (A^TS^TSA)^{-1}$, and the James-Stein estimator in (\ref{eq:shrink_Sigma}), we get the following estimator for $x_{LS}$:
\begin{equation}
    \label{eq:x_hat_s}
    \widehat{x}_S = \left( 1- \frac{(d-2) \norm{y^\perp}_2^2}{m\norm{SA\widehat{x}}_2^2} \right)\widehat{x}
\end{equation}
where $y^\perp=y-Ax_{LS}$ is as defined in Lemma \ref{lem:ls} and $\widehat{x}$ is the solution of the classical Gaussian sketch (\ref{eq:sketchls}).

However, $\norm{y^\perp}^2=\norm{y-Ax_{LS}}_2^2$ cannot be computed directly without knowing the true solution $x_{LS}$ itself. So we approximate $\norm{y^\perp}_2^2$ as follows (an alternate approximation is shown in Section \ref{sec:empirical}).
\begin{equation}
    \label{eq:yperp_approx}
    \norm{y^\perp}_2^2 \approx \frac{m-d-1}{m-1} \norm{A\widehat{x}-y}_2^2
\end{equation}

With this approximation, we use the following shrinkage estimator for $x_{LS}$:
\begin{equation}
    \label{eq:x_hat_shr}
    \widehat{x}_{shr} = \left( 1- \frac{(d-2)(m-d-1) \norm{A\widehat x-y}_2^2}{m(m-1)\norm{SA\widehat{x}}_2^2} \right)\widehat{x}
\end{equation}

\subsection{Upper Bound on Error}
\begin{theorem}
\label{theo:upper}
\begin{enumerate}[(a)]
\item The approximation in (\ref{eq:yperp_approx}) is an unbiased estimate, i.e.
\begin{equation}
    \E\left[\frac{m-d-1}{m-1} \norm{A\widehat{x}-y}_2^2\right] = \norm{y^\perp}_2^2
\end{equation}
\item For $m>d+3$, the error of the shrinkage estimator (\ref{eq:x_hat_shr}) in the $\norm{SA(\cdot)}_2$ norm is upper bounded as
\begin{equation}
    \label{eq:upper_bound_SA}
    \E\left[ \norm{SA(\widehat{x}_{shr} - x_{LS})}_2^2 \right] \leq \frac{d}{m} \norm{y^\perp}_2^2 \left( 1- \frac{1 - \epsilon'(d,m)} {1 + ({m}/{d}) \norm{Ax_{LS}}_2^2 /\norm{y^\perp}_2^2} \right)
\end{equation}
where $\epsilon'(d,m) = \frac{4(d-1)}{d^2} + \frac{2(d-2)^2}{d(m-1)(m-d-3)}$.
\item If $m>\mathcal{O}\left( (d+\log\left(1/\delta\right))/\epsilon^2\right)$, then with probability greater than $1-\delta$, the $\norm{A(\cdot)}_2$ norm of the error (prediction error) is upper bounded as
\begin{equation}
    \label{eq:upper_bound}
    \E\left[ \norm{A(\widehat{x}_{shr} - x_{LS})}_2^2 \right] \leq (1+\epsilon)\frac{d}{m} \norm{y^\perp}_2^2 \left( 1- \frac{1 - \epsilon'(d,m)} {1 + ({m}/{d}) \norm{Ax_{LS}}_2^2 /\norm{y^\perp}_2^2} \right)
\end{equation}
\end{enumerate}
\end{theorem}
Part (c) follows from part (b) since the Gaussian sketching matrix is an Oblivious Subspace Embedding \cite{Nelson2016}. 

To examine this upper bound in more detail, we look at the ratio of the prediction errors due to the shrinkage estimator and the classical Gaussian sketch.
\begin{equation}
    \frac{\E\left[\norm{A(\widehat{x}_{shr} - x_{LS})}_2^2 \right]} {\E\left[\norm{A(\widehat{x}-x_{LS})}_2^2 \right]} \leq (1+\epsilon)\frac{m-d-1}{m} \left( 1- \frac{1 - \epsilon'(d,m)} {1 + ({m}/{d}) \norm{Ax_{LS}}_2^2 /\norm{y^\perp}_2^2} \right) = R(A,y)
\end{equation}
Clearly, $R(A,y) < (1+\epsilon)$. But for large enough $m$, we expect to have $R(A,y)<1$ for any $A,y$. (Through empirical results we will see that this happens for $m$ that is not much larger than $d$). As $m$ becomes very large, $R(A,y)\to1$, indicating that the error of the shrinkage estimator approaches the classical Gaussian sketch as $m$ becomes large.

The factor $\rho(A,y)=\norm{Ax_{LS}}_2^2/ \norm{y^\perp}_2^2 =\norm{Ax_{LS}}_2^2/\norm{Ax_{LS}-y}_2^2$ can be interpreted as a signal-to-noise ratio (SNR). 
$R(A,y)$ is smaller when $\rho(A,y)$ is small, i.e. the original least squares problem is more noisy. The shrinkage factor in (\ref{eq:x_hat_shr}) also depends on the factor $\norm{SA\widehat x}_2^2/\norm{A\widehat x-y}_2^2$ which can be regarded as an approximation of $\rho(A,y)$. When this ratio is large (with respect to $d/m$), the shrinkage estimator will be very close to the classical Gaussian sketch.

Note that this upper bound can be smaller than the lower bound (\ref{eq:lower_infinity}) because it is a worst case lower bound, and the worst case occurs when $\rho(A,y)\to\infty$. Indeed there cannot be any estimator that beats the lower bound in (\ref{eq:lower_infinity}) for all values of $A,y$.

\subsection{Matching the Upper and Lower Bounds}
The lower bound in Theorem \ref{theo:lower} and the upper bound in Theorem \ref{theo:upper} have very similar forms and we show that the shrinkage estimator achieves this lower bound under certain conditions.

Suppose we know that the true solution $x_{LS}$ satisfies the SNR condition $\rho(A,y)\leq\eta^2$. If $\norm{x_{LS}}_\infty^2 \leq \eta^2\norm{y^\perp}_2^2/d\sigma_{max}(A^TA)$,
\begin{equation*}
    \rho(A,y) = 
    \frac{\norm{Ax_{LS}}_2^2}{\norm{y^\perp}_2^2} \leq \frac{\sigma_{max}(A^TA)\norm{x_{LS}}_2^2}{\norm{y^\perp}_2^2} \leq \frac{d\sigma_{max}(A^TA)\norm{x_{LS}}_\infty^2}{\norm{y^\perp}_2^2} \leq \eta^2
\end{equation*}
So we can choose $B^2=\eta^2\norm{y^\perp}_2^2/d\sigma_{max}(A^TA)$ in the lower bound (\ref{eq:lower_B}). The lower bound becomes
\begin{equation*}
    \E\left[\norm{A(\widetilde{x}-x_{LS})}_2^2 \right] \geq \frac{d}{m} \norm{y^\perp}_2^2 \left(1- \frac{d\pi^2\sigma_{max}(A^TA)} {m\eta^2\sigma_{min}(A^TA)} \right)
\end{equation*}

If $d\gg1$ and $m\gg d$ , we ignore $\epsilon'(d,m)$ and the 1 in the denominator and the upper bound (\ref{eq:upper_bound}) now becomes:
\begin{equation*}
    \E\left[ \norm{A(\widehat{x}_{shr} - x_{LS})}_2^2 \right] \leq (1+\epsilon)\frac{d}{m} \norm{y^\perp}_2^2 \left( 1- \frac{d} {m \eta^2} \right)
    \label{eq:tight_upper_bound}
\end{equation*}

Comparing this with the lower bound in (\ref{eq:lower_B}), it is achieved for large $m$ upto the constant of $\pi^2$ in the lower bound, when:
\begin{enumerate}
\item we know that $\rho(A,y)\leq\eta^2$ (or the corresponding upper bound on the $\ell^2$ or $\ell^\infty$ norm),
\item $d$ and $m$ are large enough so that $d/m$ and $\epsilon$ are small, and \item the data matrix $A$ is well-conditioned so that $\sigma_{max}(A^TA)\approx\sigma_{min}(A^TA)$.
\end{enumerate}

Condition 1 can be seen as a regularity condition on the solution which is typically desired. As will be seen in the empirical results, Condition 2 can be satisfied for $m$ that is not too large. Thus the only data-dependent condition is Condition 3.

\section{Empirical Results}
\label{sec:empirical}
The decrease in prediction error due to the shrinkage estimator will be demonstrated empirically on both simulated and real datasets. The behaviour of the shrinkage estimator is examined for different SNR values of the data, and also compared with the upper and lower bounds. In all the following plots, the (normalized) prediction error refers to $\frac{1}{n}\norm{A(\widehat{x}-x_{LS})}_2^2$ (unless otherwise specified).

\subsection{Results on Simulated Data}
\begin{figure}[!h]
    \centering
    \begin{subfigure}[b]{0.49\columnwidth}
        \centering
        \includegraphics[width=\textwidth]{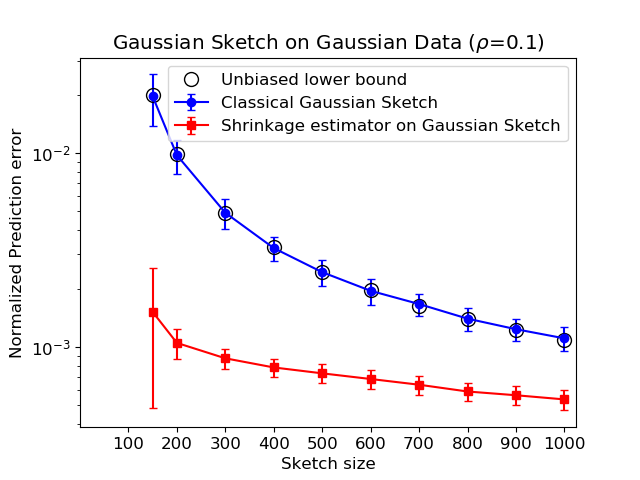}
        \caption{}
        \label{fig:gaussian_data_0.1}
    \end{subfigure}
    \begin{subfigure}[b]{0.49\columnwidth}
        \centering
        \includegraphics[width=\textwidth]{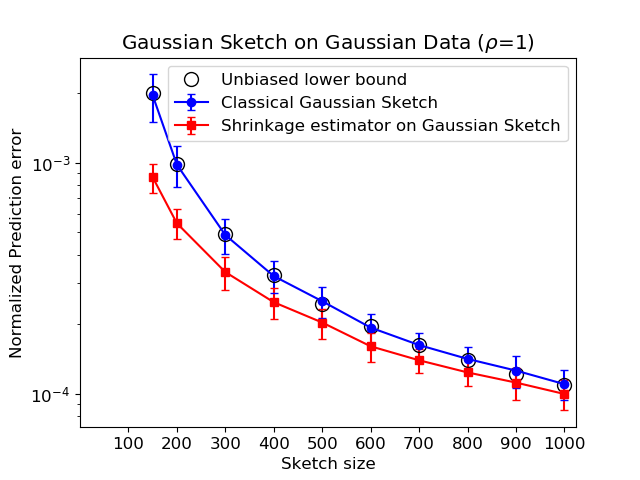}
        \caption{}
        \label{fig:gaussian_data_1}
    \end{subfigure}
    \begin{subfigure}[b]{0.49\columnwidth}
        \centering
        \includegraphics[width=\textwidth]{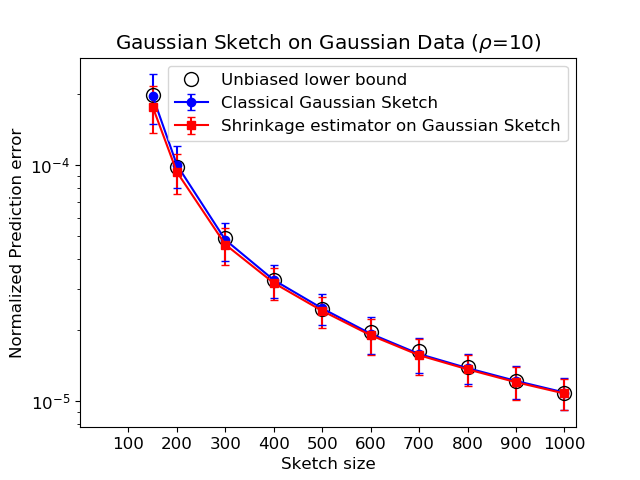}
        \caption{}
        \label{fig:gaussian_data_10}
    \end{subfigure}
    \begin{subfigure}[b]{0.49\columnwidth}
        \centering
        \includegraphics[width=\textwidth]{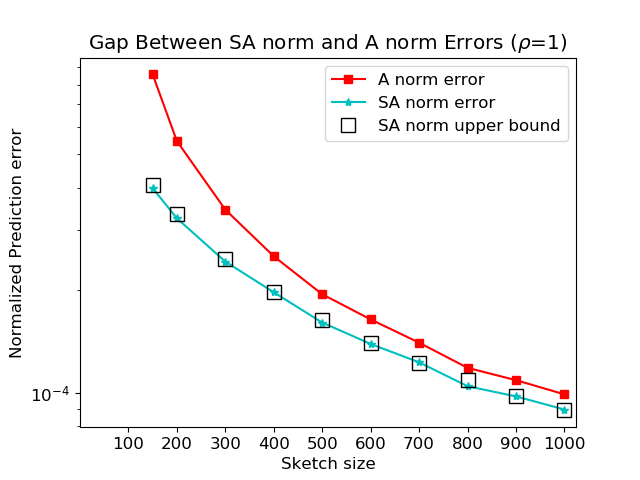}
        \caption{}
        \label{fig:gaussian_data_sa}
    \end{subfigure}
    \caption{Shrinkage estimator applied to Gaussian data. (a-c): Error $\frac{1}{n}\protect\norm{A(\widehat{x}-x_{LS})}_2^2$ for the classical Gaussian sketch as well as the shrinkage estimator, versus the sketch size ($m$), for different values of $\rho(A,y)$. (d): Comparison of the errors $\frac{1}{n}\protect\norm{A(\widehat{x}_{shr}-x_{LS})}_2^2$ and $\frac{1}{n}\protect\norm{SA(\widehat{x}_{shr}-x_{LS})}_2^2$, and the upper bound (\ref{eq:upper_bound_SA}). All the plots are averages over 100 independent sketches and the error bars show the standard deviation on each side.}
    \label{fig:gaussian_data}
\end{figure}

Figure \ref{fig:gaussian_data} compares the prediction error for the classical Gaussian sketch and the shrinkage estimator on simulated data. For these experiments, we generate each row of $A$ from a multivariate Gaussian distribution with mean $\mu_i=1$ and covariance $\Sigma_{ij}=(0.5)^{|i-j|}$. We choose $x_{LS}$ as an i.i.d. random normal vector and normalize such that $\norm{Ax_{LS}}_2=1$. Then we choose $y=Ax_{LS}+\alpha Vw$ where $V$ is a basis for the null space of $A^T$ and $w\sim\N(0,\sigma^2I_{n-d})$ (Gaussian noise model) with $n=1024$ and $d=100$. By choosing $\alpha=(\sqrt{\rho}\norm{Vw}_2)^{-1}$, we set the SNR  $\rho(A,y)=\rho$ for different values of $\rho$. This will be referred to as \textbf{Gaussian data}.

In Figure \ref{fig:gaussian_data}, the shrinkage estimator achieves much lower error, and the decrease is most noticeable at lower sketch sizes. The decrease in error due to the shrinkage estimator is not much when the SNR is large (Fig. \ref{fig:gaussian_data_10}), but it is considerable when the SNR is small (Fig. \ref{fig:gaussian_data_0.1}). These plots also show the lower bound for unbiased estimators (\ref{eq:unbiased_lower}) and that the classical Gaussian sketch achieves the lower bound. Figure \ref{fig:gaussian_data_sa} compares the errors $\frac{1}{n}\norm{A(\widehat{x}-x_{LS})}_2^2$ and $\frac{1}{n}\norm{SA(\widehat{x}-x_{LS})}_2^2$ for the shrinkage estimator, showing that the $(1+\epsilon)$ factor in (\ref{eq:upper_bound}) becomes close to 1 for a sketching dimension that is not too large. The $\frac{1}{n}\norm{SA(\widehat{x}-x_{LS})}_2^2$ error is also compared with the upper bound in (\ref{eq:upper_bound_SA}), showing that the upper bound matches the empirical error.

\subsection{Shrinkage Estimator on other Sketching Methods}
\begin{figure}[h!]
    \centering
    \begin{subfigure}[b]{0.49\columnwidth}
        \centering\includegraphics[width=\textwidth]{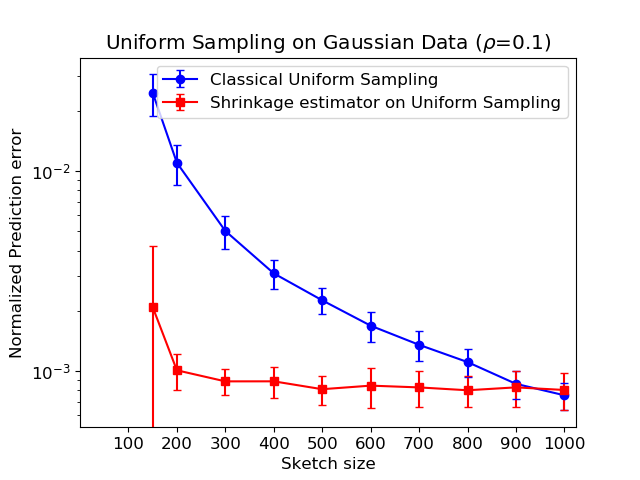}
        \caption{}
    \end{subfigure}
    \begin{subfigure}[b]{0.49\columnwidth}
        \centering\includegraphics[width=\textwidth]{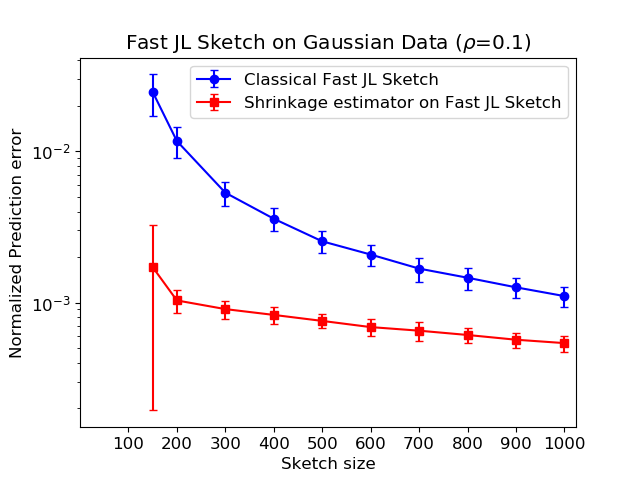}
        \caption{}
    \end{subfigure}
    \caption{Comparing the classical uniform sampling and Fast JL Transform based sketches with the shrinkage estimator applied on them. All the plots are averages over 100 independent sketches and the error bars show the standard deviation on each side.}
    \label{fig:fjlt_sketch}
\end{figure}

Although the analysis for the shrinkage estimator has been done only for the Gaussian sketch, it is seen to work equally well when applied to the estimate obtained from other types of sketching matrices such as uniform row sampling and the Fast JL transform (Figure \ref{fig:fjlt_sketch}). In fact, it also works equally well for the CW Sketch, Rademacher sketch and sampling based on row norms or leverage scores (results not shown here).

\subsection{Alternate Approximation to (\ref{eq:yperp_approx})}
Instead of approximating $\norm{y^\perp}_2^2$ using (\ref{eq:yperp_approx}), we can also use the following approximation:
\begin{equation}
    \norm{y^\perp}_2^2 \approx \frac{m}{m-d} \norm{SA\widehat x-Sy}_2^2
    \label{eq:yperp_approx_alt}
\end{equation}
This is also an unbiased estimator (see proof of Theorem \ref{theo:upper}), and is useful when only the sketched data $SA$ and $Sy$ can be be observed but not the original data $A$ and $y$. This can be the case when the measurement system itself gives a compressed measurement, or in scenarios with memory constraints. This gives us the alternate estimator
\begin{equation}
    \label{eq:x_hat_alt}
    \widehat{x}_{alt} = \left( 1- \frac{(d-2) \norm{SA\widehat x-Sy}_2^2}{(m-d)\norm{SA\widehat{x}}_2^2} \right)\widehat{x}
\end{equation}
Figure \ref{fig:alt_approx} shows that the true James-Stein estimator (\ref{eq:x_hat_s}), the original approximation (\ref{eq:x_hat_shr}) and the alternate approximation (\ref{eq:x_hat_alt}) achieve similar expected errors.

\subsection{Positive Part Shrinkage Estimator}
\begin{figure}[h!]
    \centering
    \begin{subfigure}[b]{0.49\columnwidth}
        \centering\includegraphics[width=\textwidth]{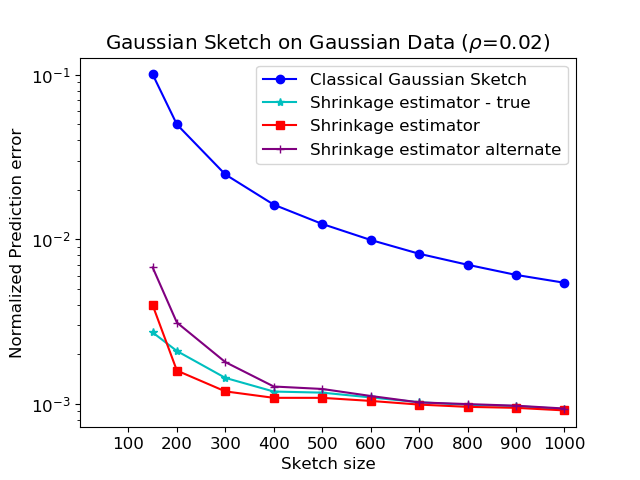}
        \caption{}
        \label{fig:alt_approx}
    \end{subfigure}
    \begin{subfigure}[b]{0.49\columnwidth}
        \centering\includegraphics[width=\textwidth]{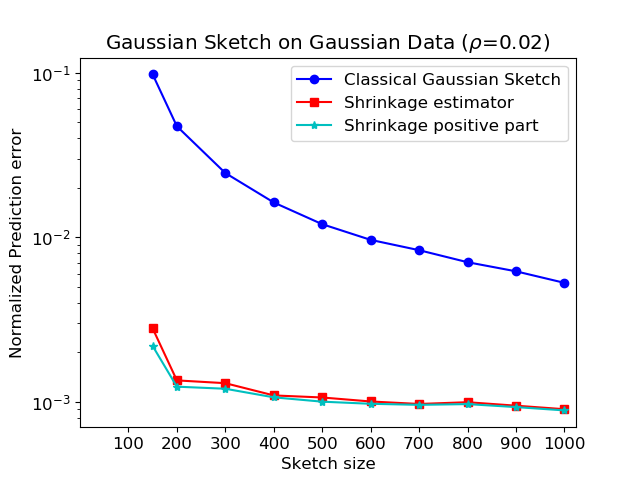}
        \caption{}
        \label{fig:pp}
    \end{subfigure}
    \caption{(a): Comparison of the true James-Stein estimator, the proposed shrinkage estimator and the alternate approximation. (b): Comparison of the shrinkage estimator and the positive part estimator. All the plots are averages over 100 independent sketches}
\end{figure}

The James-Stein estimator shrinks the estimate of the mean towards 0. However it is possible that the shrinkage factor can be too large causing the estimator to switch signs. To avoid this, we can use a positive-part James-Stein estimator which is known to have strictly lower expected squared error than the James-Stein estimator too.
\begin{equation}
    \label{eq:x_hat_pp}
    \widehat{x}_{pp} = \max\left\{\left( 1- \frac{(d-2)(m-d-1) \norm{A\widehat x-y}_2^2}{m(m-1)\norm{SA\widehat{x}}_2^2} \right),0\right\}\widehat{x}
\end{equation}
In these experiments, it was seen that there is an improvement in the error only at very low values of $\rho(A,y)$ (Figure \ref{fig:pp}).

\subsection{Results on Real Data}
\begin{figure}[h!]
    \centering
    \begin{subfigure}[b]{0.49\columnwidth}
        \centering
        \includegraphics[width=\textwidth]{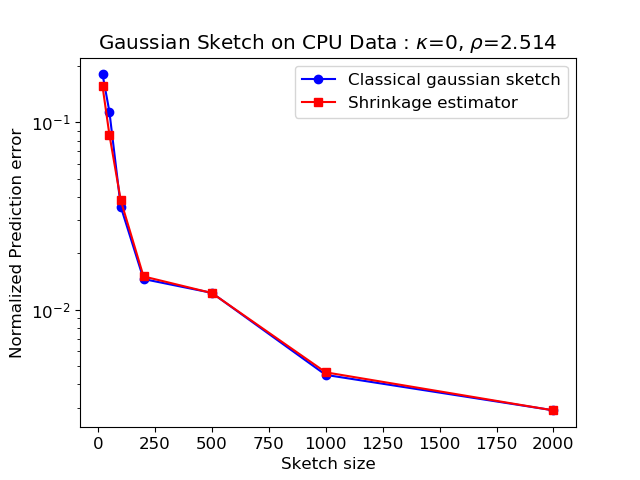}
        \caption{}
    \end{subfigure}
    \begin{subfigure}[b]{0.49\columnwidth}
        \centering
        \includegraphics[width=\textwidth]{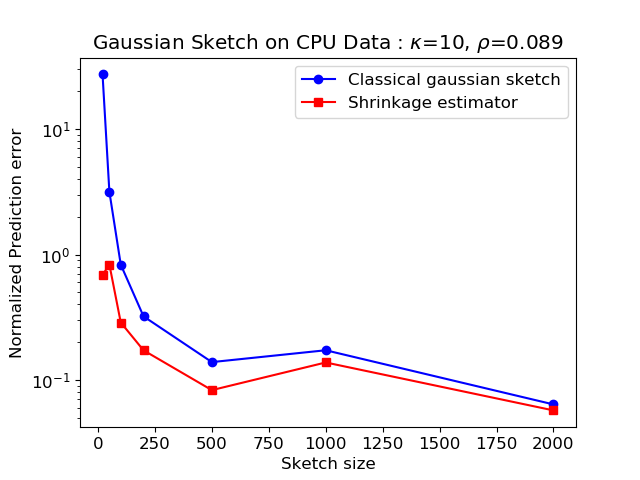}
        \caption{}
    \end{subfigure}
    \begin{subfigure}[b]{0.49\columnwidth}
        \centering
        \includegraphics[width=\textwidth]{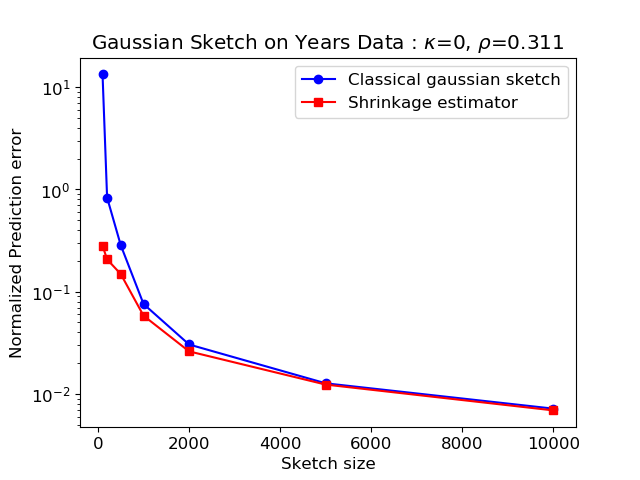}
        \caption{}
    \end{subfigure}
    \begin{subfigure}[b]{0.49\columnwidth}
        \centering
        \includegraphics[width=\textwidth]{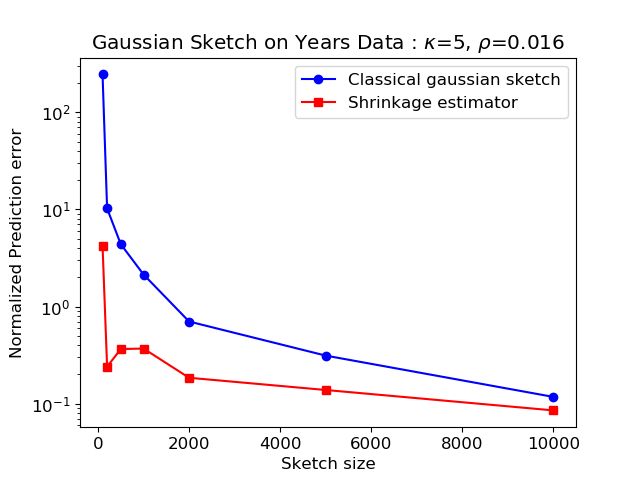}
        \caption{}
    \end{subfigure}
    \caption{Results of shrinkage estimator and classical Gaussian sketch on the CPU data ($n=8192,d=12$) and Years data ($n=463715,d=90$). $\kappa$ is the variance of the added noise scaled by $\protect\norm{y^\perp}_2^2$ and $\rho$ is the SNR.}
    \label{fig:cpuyears}
\end{figure}

We present empirical results on 4 real datasets. The CPU dataset ($n=8192,d=12$) and the Years dataset ($n=463715,d=90$) \cite{libsvm} are regression datasets. We compare the shrinkage estimator and classical Gaussian sketch for both these datasets with varying levels of added noise, i.e. we add Gaussian noise of variance $\kappa\norm{y^\perp}_2^2$ to $y$ to vary the SNR $\rho$. The results are seen in Figure \ref{fig:cpuyears}a,b for CPU data and in Figure \ref{fig:cpuyears}c,d for the Years data. The shrinkage estimator performs much better than the classical sketch when the data is very noisy and when $n\gg d$, and performs equally well as the classical sketch otherwise.

We also show results for least-squares based classification on the MNIST and SVHN datasets. For SVHN ($n=20000,d=3072$), we classify between two digits using vectorized images in the rows of $A$ and 0-or-1 labels in $y$. For MNIST ($n=48000,d=1000$), we perform a multiclass classification with vectorized images passed through a randomly initialized ReLU layer in the rows of $A\in\R^{n\times d}$, and $Y\in\R^{n\times10}$ is a matrix with one-hot vector labels. Thus we have a Frobenius norm regression problem: minimize $\norm{AX-Y}_F^2$. With $\widehat X$ as the solution of the sketched version (minimize $\norm{SAX-SY}_F^2$), the shrinkage estimator in this case is simply 
\begin{equation}
    \label{eq:x_hat_shr_fro}
    \widehat{X}_{shr} = \left( 1- \frac{(d-2)(m-d-1) \norm{A\widehat X-Y}_F^2}{m(m-1)\norm{SA\widehat{X}}_F^2} \right)\widehat{X}
\end{equation}
This estimator satisfies the same error upper bound (\ref{eq:upper_bound}) with the vector norms replaced by the Frobenius norms.

\sloppy

\begin{figure}[h!]
    \centering
    \begin{subfigure}[b]{0.49\columnwidth}
        \centering
        \includegraphics[width=\textwidth]{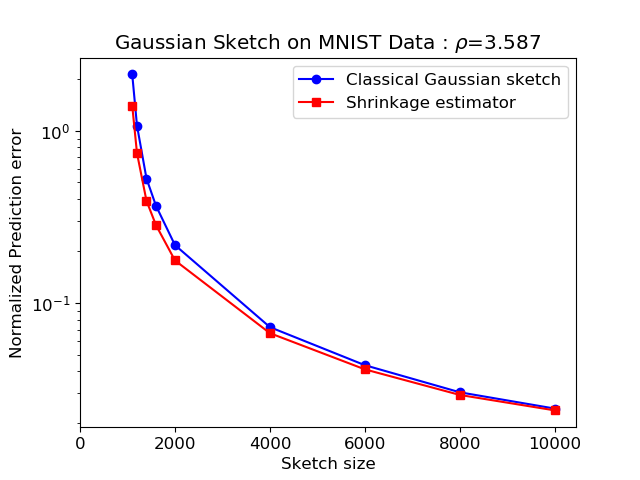}
        \caption{}
    \end{subfigure}
    \begin{subfigure}[b]{0.49\columnwidth}
        \centering
        \includegraphics[width=\textwidth]{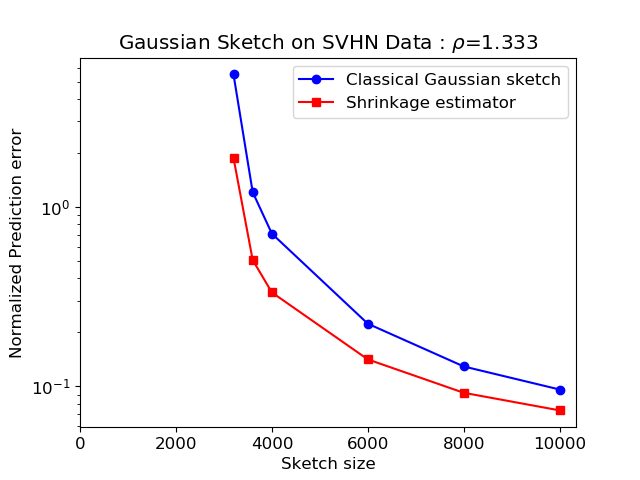}
        \caption{}
    \end{subfigure}
    \caption{Results of shrinkage estimator and classical Gaussian sketch on the MNIST ($n=48000,d=1000$) and SVHN data ($n=20000,d=3072$). $\rho$ is the SNR.}
    \label{fig:mnistsvhn}
\end{figure}

Figure \ref{fig:mnistsvhn} shows the prediction error on these two datasets. Since $\norm{A\widehat x-y}_2^2 = \norm{A(\widehat x-x_{LS})}_2^2 + \norm{Ax_{LS}-y}_2^2$, the mean squared error loss is only a constant added to the prediction error. While the SNR of MNIST is high and the shrinkage estimator performs only slightly better than the classical Gaussian sketch, the shrinkage estimator performs much better than the classical Gaussian sketch in SVHN.

\fussy

\section{Conclusion and Further Directions}
\label{sec:conclusion}
Applying the James-Stein estimator to sketched least squares optimization results in smaller prediction error, especially when the data is more noisy. The estimator is easy to compute based on the classical sketching methods, and can allow us to reduce the sketch sizes greatly. The analysis is interesting as it looks at an optimization problem as a statistical estimation problem. This idea works irrespective of the distribution of the data. Empirically it works for many common sketching methods as well. From prior results that show that common sketching solutions asymptotically converge to the Gaussian sketch solution in distribution \cite{ahfock2017statistical, li2006very}, we can argue that our results apply to those sketches as well.

The lower bounds derived in this project are new, and easy to compute, unlike existing ones such as in \cite{pilanci2014iterative}. We proved that the classical sketching method is the optimum unbiased estimator based on the Gaussian sketch. Additionally our lower bound can incorporate prior knowledge about the solution through constraints. The error of the shrinkage estimator can be made arbitrarily close to the lower bound for large $m$ and $d$, when the data is well conditioned (upto the factor of $\pi^2$), so the bounds are tight.

Since the shrinkage estimator shrinks the sketched solution towards 0, it can be regarded as a regularizer. To the best of our knowledge, there are no results showing whether an $\ell^2$ norm regularized sketching solution strictly dominates the classical sketching solution from a statistical perspective. But our shrinkage estimator is shown to strictly dominate the classical Gaussian sketch solution. It would be interesting to study the generalization properties of the shrinkage estimator proposed here, and compare it with standard regularization methods.

Another interesting direction is to apply such shrinkage estimators for right sketching methods for the minimum norm solution of underdetermined equations. Since the minimum norm solution also has the same form and has a similar distribution under Gaussian sketch, such an estimator would give lower error in that case as well.

\bibliography{bibliography/main}
\bibliographystyle{IEEEtran}

\section{Proofs}
\label{sec:proofs}
\subsection{Proof of Lemma \ref{lem:ls}}
The solution $x_{LS}$ of (\ref{eq:ls}) satisfies $\nabla \norm{Ax-y}_2^2=0$, i.e. $A^T(Ax_{LS}-y)=0$. By solving this equation we get $x_{LS}=(A^TA)^{-1}A^Ty=A^\dagger y$ (assuming that $A$ is full column rank, so that $A^TA$ is invertible). By setting $y^\perp=y-Ax_{LS}$, we get $A^Ty^\perp=0$. Similarly, the sketched solution $\widehat{x}$ of (\ref{eq:sketchls}) satisfies $(SA)^T(SA\widehat{x}-Sy)=0$. Solving this, we get
\begin{align*}
    \widehat{x} &= (A^TS^TSA)^{-1}A^TS^TSy \\
    &= (A^TS^TSA)^{-1}A^TS^TS(Ax_{LS} + y^\perp) \\
    &= x_{LS} + (A^TS^TSA)^{-1}A^TS^TSy^\perp \\
    &= x_{LS} + (SA)^\dagger Sy^\perp
\end{align*}
Since each element of $SA$ and $Sy^\perp$ is a sum of independent Gaussians, each column of $SA$, and $Sy^\perp$ follow a multivariate Gaussian distribution. In particular, $Sy^\perp\sim\N(0, \frac{1}{m}\norm{y^\perp}_2^2 I_m)$. Further, \begin{equation*}
    \E[(SA)^TSy^\perp] = \E[A^TS^TSy^\perp] = A^T\E[S^TS]y^\perp = A^Ty^\perp = 0 \quad \text{since } \E[S^TS] = I
\end{equation*}
Thus, $SA$ and $Sy^\perp$ are uncorrelated Gaussian random variables, so they are independent too. Finally, we note that conditioned on $SA$, $\widehat{x}$ must have a Gaussian distribution (since $Sy^\perp$ is Gaussian) with $\E[\widehat{x}\mid SA] = x_{LS}$ since $\E[Sy^\perp\mid SA]=\E[Sy^\perp]=0$, and
\begin{align*}
\text{Cov}[\widehat{x}\mid SA] &=
\E[(\widehat{x}-x_{LS})(\widehat{x}-x_{LS})^T \mid SA] \\
&= (SA)^\dagger \E[(Sy^\perp)(Sy^\perp)^T \mid SA] (SA)^{\dagger T}\\
&= (SA)^\dagger\text{Cov}[Sy^\perp\mid SA](SA)^{\dagger T} \\ &= (SA)^\dagger\text{Cov}[Sy^\perp](SA)^{\dagger T}\\
&= (A^TS^TSA)^{-1}A^TS^T \left(\frac{1}{m}\norm{y^\perp}_2^2 I_m\right) SA(A^TS^TSA)^{-1}\\
&= \frac{1}{m}\norm{y^\perp}_2^2 (A^TS^TSA)^{-1}
\end{align*}
\endproof{}
\subsection{Proof of Lemma \ref{lem:gaussian_sketch}}
From Lemma \ref{lem:ls},
\begin{equation*}
    \widehat{x} \sim \N\left(x_{LS}, \frac{1}{m}\norm{y^\perp}_2^2 (A^TS^TSA)^{-1}\right) \mid SA
\end{equation*}
Hence we have
\begin{align*}
    \E\left[\norm{A(\widehat x-x_{LS})}_2^2 \mid SA \right] &= \E\,\text{tr} \left[A(\widehat x-x_{LS})(\widehat x-x_{LS})^TA^T \mid SA \right] \\
    &= \text{tr} \left[A \E\left[(\widehat x-x_{LS})(\widehat x-x_{LS})^T\mid SA \right]A^T\right]\\
    &= \text{tr} \left[ A\left(\frac{1}{m}\norm{y^\perp}_2^2 (A^TS^TSA)^{-1}\right)A^T\right]
\end{align*}
Taking expectation with respect to $SA$ on both sides,
\begin{align*}
    \E\left[\norm{A(\widehat x-x_{LS})}_2^2 \right] &= \E_{SA}\;\text{tr}\left[ A\left(\frac{1}{m}\norm{y^\perp}_2^2 (A^TS^TSA)^{-1}\right)A^T\right]\\
    &= \frac{1}{m}\norm{y^\perp}_2^2 \text{tr}\left[A \E_{SA}\left[(A^TS^TSA)^{-1}\right] A^T\right]
\end{align*}

Since each row of $SA$ has an independent multivariate Gaussian distribution $\N\left(0,\frac{1}{m}A^TA\right)$, $A^TS^TSA$ follows a Wishart distribution $\W_d\left(m, \frac{1}{m}A^TA\right)$, and so its inverse follows an inverted Wishart distribution that satisfies $\E\left[(A^TS^TSA)^{-1}\right]=\frac{m}{m-d-1}(A^TA)^{-1}$ \cite{gupta2018matrix}.
\begin{align*}
    \E\left[\norm{A(\widehat x-x_{LS})}_2^2 \right] &= \frac{1}{m}\norm{y^\perp}_2^2 \text{tr}\left[A \frac{m}{m-d-1}(A^TA)^{-1} A^T\right]\\
    &= \frac{1}{m-d-1} \norm{y^\perp}_2^2 \text{tr}\left[A(A^TA)^{-1}A^T\right]\\
    &= \frac{1}{m-d-1} \norm{y^\perp}_2^2 \text{tr}\left[(A^TA)^{-1}A^TA\right]\\
    &= \frac{d}{m-d-1} \norm{y^\perp}_2^2 
\end{align*}
\endproof{}

\subsection{Proof of Lemma \ref{lem:stein}}
The proof is a consequence of Stein's lemma from \cite{stein1981} which can be proved using the integration by parts formula for multivariable calculus. With $\widehat \theta$ as in (\ref{eq:shrink_Sigma}),
\begin{align*}
    \E[(\widehat\theta - \theta)^T\Sigma^{-1} (\widehat\theta - \theta)] &=
    \E\left[ \left(X-\theta - \frac{(d-2)X}{X^T\Sigma^{-1}X}\right)^T \Sigma^{-1} \left(X-\theta - \frac{(d-2)X}{X^T\Sigma^{-1}X}\right) \right] \\
    &= \E[(X-\theta)^T\Sigma^{-1}(X-\theta)] - 2(d-2) \E\left[\frac{(X-\theta)^T\Sigma^{-1}X}{X^T\Sigma^{-1}X} \right]\\
    & \qquad + (d-2)^2 \E\left[ \frac{1}{X^T\Sigma^{-1}X} \right]
\end{align*}
The first term can be simplified as
\begin{align*}
    \E[(X-\theta)^T\Sigma^{-1}(X-\theta)] &= 
    \E\,\text{tr}[(X-\theta)^T\Sigma^{-1}(X-\theta)] \\
    &= \E\,\text{tr}[\Sigma^{-1}(X-\theta)(X-\theta)^T] \\
    &= \text{tr}(\Sigma^{-1}\E[(X-\theta)(X-\theta)^T])\\
    &= \text{tr}(\Sigma^{-1}\Sigma)\\ &= d
\end{align*}
For the second term, we use the following identity from multivariable calculus \cite{rogers2011calculus}:
\begin{align*}
    \int_\Omega u\,\nabla\cdot\mathbf V\,d\Omega  \ =\ \int_{\Gamma} u \mathbf V \cdot \hat{\mathbf n}\,d\Gamma - \int_\Omega  \nabla(u)\cdot\mathbf V\,d\Omega
\end{align*}
where $\Omega\subseteq\R^d$, $u:\Omega\to\R$, $\mathbf V:\Omega\to\R^d$, $\Gamma$ is the boundary of $\Omega$ and $\hat{\mathbf n}$ is the outward normal to $\Gamma$. $\nabla(u)$ is the gradient of $u$ and $\nabla\cdot\mathbf V$ is the divergence of $\mathbf V$.

We choose $\Omega=\R^d$, $u=f(x)$ as the density function of $\N(\theta,\Sigma)$ such that $\nabla(u)=-f(x)\Sigma^{-1}(x-\theta)$, and $V=\frac{x}{x^T\Sigma^{-1}x}$. Since the density function goes to zero at the boundary of $\Omega$, the first integral on the right side is zero.
\begin{align*}
    \E\left[\frac{(X-\theta)^T\Sigma^{-1}X}{X^T\Sigma^{-1}X} \right] &= \int_\Omega\,f(x) \frac{(x-\theta)^T\Sigma^{-1}x}{x^T\Sigma^{-1}x}\,d\Omega\\
    &= \int_\Omega -\nabla f(x)\cdot \frac{x}{x^T\Sigma^{-1}x} \,d\Omega \\
    &= \int_\Omega f(x)\,\nabla\cdot\left( \frac{x}{x^T\Sigma^{-1}x}\right)\,d\Omega\\
    &= \E\left[ \nabla\cdot\left( \frac{X}{X^T\Sigma^{-1}X}\right) \right]\\
    &= \E\left[ \sum_{i=1}^d \frac{X^T\Sigma^{-1}X - 2X_i(\Sigma^{-1}X)_i}{(X^T\Sigma^{-1}X)^2} \right]\\
    &= (d-2)\E\left[ \frac{1}{X^T\Sigma^{-1}X} \right]
\end{align*}

Putting these parts together, we get
\begin{align*}
    \E[(\widehat\theta - \theta)^T\Sigma^{-1} (\widehat\theta - \theta)] &=
    d - (d-2)^2\E\left[ \frac{1}{X^T\Sigma^{-1}X} \right]
\end{align*}
\endproof{}

\subsection{Proof of Theorem \ref{theo:unbiased_lower}}
We first note by the definition of $y^\perp$ that \[ \widetilde y = Sy = SAx_{LS} + Sy^\perp \] where $SAx_{LS}$ and $Sy^\perp$ are both Gaussian random vectors, and are independent to each other (from Lemma \ref{lem:ls}). So, conditioned on $SA$, $\widetilde y$ is also a Gaussian random vector with mean $SAx_{LS}$ and the vector $Sy^\perp$ acting as independent, zero-mean Gaussian noise, i.e.
\begin{equation}
    \widetilde y \sim \N\left( SAx_{LS}, \frac{1}{m}\norm{y^\perp}_2^2 I_m \right) \mid SA
\end{equation}

Our problem is to estimate the parameter $x_{LS}$ in the mean of this Gaussian distribution. Let $f(\widetilde y;x_{LS})$ be the density function of the Gaussian shown above, parameterized by $x_{LS}$. Then we define the Fisher information matrix for estimation of $x_{LS}$ from $\widetilde y$ as
\begin{align*}
    I(\widetilde y;x_{LS}) &= \E \left[ \nabla_{x_{LS}} \log f(\widetilde y;x_{LS}) \nabla_{x_{LS}} \log f(\widetilde y;x_{LS})^T \right]
\end{align*}
Since $\log f(\widetilde y;x_{LS}) = -(1/2) (\widetilde y-SAx_{LS})^T \left(\frac{1}{m} \norm{y^\perp}_2^2 I_m \right)^{-1} (\widetilde y-SAx_{LS}) + $ constant,
\begin{align*}
    I(\widetilde y;x_{LS}) &= \E\left[ (SA)^T \left(\frac{1}{m} \norm{y^\perp}_2^2 I_m \right)^{-1} (\widetilde y-SAx_{LS}) (\widetilde y-SAx_{LS})^T \left(\frac{1}{m} \norm{y^\perp}_2^2 I_m \right)^{-1} SA \right]\\
    &= \frac{m}{\norm{y^\perp}_2^2} A^TS^TSA \quad \text{since } \E\left[(\widetilde y-SAx_{LS}) (\widetilde y-SAx_{LS})^T\right]=\frac{1}{m} \norm{y^\perp}_2^2 I_m
\end{align*}

For unbiased estimators, we have a lower bound which is easily derived using the Cram\'er-Rao lower bound. According to the following version of the Cram\'er-Rao lower bound which only holds for unbiased estimators \cite[page~147]{borovkov},
\begin{align*}
    \E\left[(\widetilde x-x_{LS})(\widetilde x-x_{LS})^T \mid SA \right] \succeq I(\widetilde y;x_{LS})^{-1}
\end{align*}
where the inequality $A\succeq B$ means that $A-B$ is positive semi-definite. Since the Gaussian distribution and the Fisher information above are conditioned on $SA$, this lower bound holds under the same conditioning. We can lower bound the mean squared error as follows:
\begin{align*}
    \E\left[\norm{A(\widetilde x-x_{LS})}_2^2 \mid SA \right] &= \text{tr} \left[A\E\left[(\widetilde x-x_{LS})(\widetilde x-x_{LS})^T\mid SA \right] A^T\right]\\
    &\geq \frac{\norm{y^\perp}_2^2}{m} \text{tr}\left[A (A^TS^TSA)^{-1} A^T \right]
\end{align*}

After taking expectation with respect to $SA$ on both sides, we follow the same steps as in the proof of Lemma \ref{lem:gaussian_sketch} to get
\begin{align*}
    \E\left[\norm{A(\widetilde x-x_{LS})}_2^2 \right] &\geq \frac{d}{m-d-1} \norm{y^\perp}_2^2
\end{align*}
\endproof{}

\subsection{Proof of Theorem \ref{theo:lower}}
Define the Fisher information matrix $I(\widetilde y;x_{LS})$ as in the proof for Theorem \ref{theo:unbiased_lower}. In this proof, we will have to use slightly different and more sophisticated tools because the Cram\'er-Rao bound used in the proof for Theorem \ref{theo:unbiased_lower} holds only for unbiased estimators.

Let $\lambda(x_{LS})$ be the density function of some prior distribution for $x_{LS}$, which has support on $[-B,B]^d$ and let $I(\lambda)$ be the matrix defined by
\[ I(\lambda) = \E\left[ \nabla\log\lambda(x_{LS}) \nabla\log\lambda(x_{LS})^T \right] \]

Let $\widetilde x$ be an estimator of $x_{LS}$ derived from $\widetilde y$. Let $\Lambda=(A^TA)^{-1}$. Our problem is to estimate the parameter $x_{LS}\in\C$ such that $[-B,B]^d\subseteq\C$. We use the following special case of Theorem 1 from \cite{gill1995}, which holds under certain regularity conditions which are explained therein. These regularity conditions are easily satisfied by the Gaussian distribution. \cite{gill1995} states a lower bound on the expected error averaged with respect to the density $\lambda(x_{LS})$. Since we do not consider any probabilistic model on the data, we can lower bound the worst case error by the average error.
\begin{align*}
    \E\left[(\widetilde x-x_{LS})^T\Lambda^{-1} (\widetilde x-x_{LS}) \mid SA \right] \geq \frac{d^2} {\E_{x_{LS}}\text{tr}\left[ \Lambda^TI(\widetilde y;x_{LS})\right] + \E_{x_{LS}}\text{tr}\left[\Lambda^TI(\lambda)\right]}
\end{align*}

In this case, the two quantities in the denominator do not depend on $x_{LS}$ and hence the expectations with respect to $x_{LS}$ are trivial. We then take the expectation of both sides with respect to $SA$.
\begin{align*}
\E\left[\norm{A(\widetilde x - x_{LS})}_2^2 \right] &= \E_{SA}\left[ \E\left[\norm{A(\widetilde x - x_{LS})}_2^2 \mid SA \right] \right] \\
&\geq \E_{SA}\left[ \frac{d^2} {\text{tr}\left[ \Lambda^TI(\widetilde y;x_{LS})\right] + \text{tr}\left[\Lambda^TI(\lambda)\right]} \right]\\
&\geq \frac{d^2} {\E_{SA}\,\text{tr}[ \Lambda^TI(\widetilde y;x_{LS})] + \text{tr}\left[\Lambda^TI(\lambda)\right]} & &\text{from }\E\frac{1}{X}\geq\frac{1}{\E X}
\end{align*}

We upper bound the two quantities in the denominator as follows:
\begin{align*}
    \E_{SA}\,\text{tr}[ \Lambda^TI(\widetilde y;x_{LS})] &= \frac{m}{\norm{y^\perp}_2^2} \E_{SA}\,\text{tr}[ (A^TA)^{-1} A^TS^TSA ]\\
    &= \frac{m}{\norm{y^\perp}_2^2} \text{tr} [(A^TA)^{-1} \E_{SA}[A^TS^TSA]]\\
    &= \frac{m}{\norm{y^\perp}_2^2} \text{tr} [(A^TA)^{-1} A^TA] & & \text{since }\E[S^TS]=I\\
    &= \frac{md}{\norm{y^\perp}_2^2}
\end{align*}
\begin{align*}
    \text{tr}[ \Lambda^TI(\lambda)] &= \text{tr} [(A^TA)^{-1}I(\lambda)]\\ &\leq \norm{(A^TA)^{-1}}_2 \text{tr}[I(\lambda)] & \text{because $(A^TA)^{-1}$ is symmetric and $I(\lambda)\succeq0$}\\
    &= \frac{1}{\sigma_{min}(A^TA)}\text{tr}[I(\lambda)] & \norm{(A^TA)^{-1}}_2 = \sigma_{max}(A^TA)^{-1} = \left(\sigma_{min}(A^TA)\right)^{-1}
\end{align*}
The minimum value of $\text{tr}[I(\lambda)]$ over all distributions $\lambda(.)$ which satisfy the regularity conditions and have support in $[-B,B]^d$, is $d\pi^2/B^2$. This result can be found in \cite[page~172]{borovkov} and has been used in \cite{barnes2019} as well. Since the worst case bound holds for any prior distribution, we use this minimum value.
\begin{align*}
    \E\left[\norm{A(\widetilde x - x_{LS})}_2^2 \right] &\geq \frac{d^2} {\frac{md}{\norm{y^\perp}_2^2} + \frac{d\pi^2} {B^2\sigma_{min}(A^TA)}} \\
    &\geq \frac{d}{m}\norm{y^\perp}_2^2 \left(1- \frac{\pi^2\norm{y^\perp}_2^2} {mB^2\sigma_{min}(A^TA)} \right) & &\text{from }(1+x)^{-1}\geq 1-x
\end{align*}

This bound holds for any $B>0$ such that the hypercube $[-B,B]^d$ lies within the set of feasible $x_{LS}$ values. For unconstrained optimization, we can take $B\to\infty$ and we get
\begin{align*}
    \E\left[\norm{A(\widetilde x - x_{LS})}_2^2 \right] &\geq \frac{d}{m} \norm{y^\perp}_2^2
\end{align*}
\endproof{}

\subsection{Proof of Theorem \ref{theo:upper}}
From Lemma \ref{lem:ls}, we have
\begin{equation*}
    \widehat{x} \sim \N\left(x_{LS}, \frac{1}{m}\norm{y^\perp}_2^2 (A^TS^TSA)^{-1}\right) \mid SA
\end{equation*}
We use the estimator $\widehat\theta(X)$ from Lemma \ref{lem:stein} with $X=\widehat x$, $\theta=x_{LS}$, and $\Sigma=\frac{1}{m}\norm{y^\perp}_2^2 (A^TS^TSA)^{-1}$ to get the estimator
\begin{align*}
    \widehat{x}_S &= \left( 1- \frac{(d-2)}{X^T\Sigma^{-1}X}\right)X \\
    &= \left( 1- \frac{(d-2) \norm{y^\perp}_2^2}{m\norm{SA\widehat{x}}_2^2} \right)\widehat{x}
\end{align*}
We look at the approximation of $\norm{y^\perp}_2^2$. For proving part (a), 
\begin{align*}
    \E\left[ \norm{A\widehat{x}-y}_2^2\right] &= \norm{Ax_{LS}-y}_2^2 + \E\left[ \norm{A(\widehat{x}-x_{LS})}_2^2 \right] & &\text{ since $A^Ty^\perp=0$ from Lemma \ref{lem:ls}} \\
    &= \norm{y^\perp}_2^2 + \frac{d}{m-d-1}\norm{y^\perp}_2^2 & &\text{ from Lemma \ref{lem:gaussian_sketch}}\\
    &= \frac{m-1}{m-d-1} \norm{y^\perp}_2^2
\end{align*}

Alternatively, we can also approximate $\norm{y^\perp}_2^2$ as $\frac{m}{m-d}\norm{SA\widehat{x}-Sy}_2^2$ (see Section \ref{sec:empirical}). This is also an unbiased estimator, as seen below. Firstly, $(SA)^T(Sy-SA\widehat x)=0$ from the optimality condition for the sketched problem. Then,
\begin{align*}
    \E\left[ \norm{SA\widehat{x}-Sy}_2^2\right] &= \E\left[\norm{SAx_{LS}-Sy}_2^2\right] - \E\left[ \norm{SA(\widehat{x}-x_{LS})}_2^2 \right]  \\
    &= \norm{Ax_{LS}-y}_2^2 - \E_{SA}\left[\E\left[ \norm{SA(\widehat{x}-x_{LS})}_2^2\mid SA\right] \right] \qquad \text{since } \E[S^TS]=I \\
    &= \norm{y^\perp}_2^2 - \frac{d}{m}\norm{y^\perp}_2^2 \qquad \text{since } \widehat x \sim \N\left(x_{LS}, \frac{1}{m}\norm{y^\perp}_2^2 (A^TS^TSA)^{-1}\right) \mid SA \\
    &= \frac{m-d}{m} \norm{y^\perp}_2^2
\end{align*}

For the purpose of deriving the upper bound in part (b), we will assume that we make two independent draws of the random matrix $S$ : $S_1$ and $S_2$. We will first solve the sketched problem using $S_1$ to obtain a solution $\widehat{x}_1$, and use that to compute $\norm{A\widehat{x}_1-y}_2^2$, for our approximation of $\norm{y^\perp}_2^2$. Then we solve the sketched problem using $S_2$ to obtain a solution $\widehat{x}_2$, and then compute the final estimator as
\begin{equation*}
    \widehat{x}_{shr} = \left( 1- \frac{(d-2)(m-d-1) \norm{A\widehat x_1-y}_2^2}{m(m-1)\norm{S_2A\widehat{x}_2}_2^2} \right)\widehat{x}_2
\end{equation*}

Note that this is only a trick used to simplify the derivation of the error bound. In practice, it is enough to use a single instance of $S$ to compute the estimator.

Let us start by analyzing the approximation term. Since $y^\perp=y-Ax_{LS}$ and $A^Ty^\perp=0$,
\begin{align*}
    \norm{A\widehat x_1-y}_2^2 &= \norm{y^\perp}_2^2 + \norm{A(\widehat{x}_1-x_{LS})}_2^2 \\
    \norm{A\widehat x_1-y}_2^4 &= \norm{y^\perp}_2^4 + 2 \norm{y^\perp}_2^2 \norm{A(\widehat{x}_1-x_{LS})}_2^2 + \norm{A(\widehat{x}_1-x_{LS})}_2^4 \\
    \E\left[\norm{A\widehat x_1-y}_2^4\right] &= \norm{y^\perp}_2^4 + \frac{2d}{m-d-1}\norm{y^\perp}_2^4 + \E\left[\norm{A(\widehat{x}_1-x_{LS})}_2^4\right] \quad \text{ from Lemma \ref{lem:gaussian_sketch}}
\end{align*}
\begin{align}
    \E\left[\left(\frac{m-d-1}{m-1} \norm{A\widehat x_1-y}_2^2\right)^2\right] &= \frac{(m-d-1)(m+d-1)}{(m-1)^2} \norm{y^\perp}_2^4 \nonumber \\
    & \qquad + \frac{(m-d-1)^2}{(m-1)^2} \E\left[\norm{A(\widehat{x}_1-x_{LS})}_2^4\right] \nonumber \\
    &= \norm{y^\perp}_2^4 - \frac{d^2}{(m-1)^2} \norm{y^\perp}_2^4 + \frac{(m-d-1)^2}{(m-1)^2} \E\left[\norm{A(\widehat{x}_1-x_{LS})}_2^4\right]
    \label{eq:approx_error}
\end{align}
To compute the last term on the right hand side, let $z=A(\widehat{x}_1-x_{LS})$, then \[z\sim\N\left(0, K=\frac{1}{m}\norm{y^\perp}_2^2 A(A^TS_1^TS_1A)^{-1}A^T\right) \mid S_1A\]We require the \nth{4} order moments of this Gaussian, which is as follows:
\begin{align*}
    \E\left[\norm{A(\widehat{x}_1-x_{LS})}_2^4 \mid S_1A \right] &= \E\left[\sum_{i=1}^n\sum_{j=1}^n z_i^2z_j^2 \right] \\
    &= \sum_{i=1}^n\sum_{j=1}^n \left( 2K_{ij}^2 + K_{ii}K_{jj} \right)
\end{align*}
Since each row of $S_1A$ has an independent multivariate Gaussian distribution $\N\left(0,\frac{1}{m}A^TA\right)$, $A^TS_1^TS_1A$ follows a Wishart distribution $\W_d\left(m, \frac{1}{m}A^TA\right)$, and so its inverse follows an inverted Wishart distribution $\W_d^{-1}\left(m, m(A^TA)^{-1}\right)$, and therefore \[K\sim\W_d^{-1}\left(m, \norm{y^\perp}_2^2 A(A^TA)^{-1}A^T\right)\] Let $P=A(A^TA)^{-1}A^T$ be the ``hat matrix". Using the covariance of elements of an inverted Wishart matrix from \cite{gupta2018matrix},
\begin{align*}
    \E\left[\norm{A(\widehat{x}_1-x_{LS})}_2^4\right] &= \E_{S_1A}\left[ \sum_{i=1}^n\sum_{j=1}^n \left( 2K_{ij}^2 + K_{ii}K_{jj} \right) \right]\\
    &= \sum_{i=1}^n\sum_{j=1}^n \left( 2\E\left[K_{ij}\right]^2 + 2Var\left(K_{ij}\right) + \E\left[K_{ii}\right] \E\left[K_{jj}\right] + Cov\left(K_{ii},K_{jj}\right) \right)\\
    &= \norm{y^\perp}_2^4\left( \frac{2\norm{P}_F^2}{(m-d-1)^2}
    + \frac{2(m-d+1)\norm{P}_F^2 + 2(m-d-1)\tr(P)^2}{(m-d)(m-d-1)^2(m-d-3)} \right. \\ & \qquad + \left. \frac{\tr(P)^2}{(m-d-1)^2} + \frac{2\tr(P)^2 + 2(m-d-1)\norm{P}_F^2}{(m-d)(m-d-1)^2(m-d-3)}\right) \\
    &= \frac{\norm{y^\perp}_2^4\left(2\norm{P}_F^2 + \tr(P)^2 \right)}{(m-d-1)(m-d-3)}
\end{align*}
For the hat matrix, $\tr(P)=\tr(A(A^TA)^{-1}A^T) = \tr((A^TA)^{-1}A^TA)=d$, and $\norm{P}_F^2 = \tr(P^TP) = \tr(P) = d$.
\begin{align*}
    \E\left[\norm{A(\widehat{x}_1-x_{LS})}_2^4\right] &= \frac{\norm{y^\perp}_2^4 d(d+2)}{(m-d-1)(m-d-3)}
\end{align*}
Substituting this in (\ref{eq:approx_error}),
\begin{align}
    \E\left[\left(\frac{m-d-1}{m-1} \norm{A\widehat x_1-y}_2^2\right)^2\right] &=
    \norm{y^\perp}_2^4 - \frac{d^2}{(m-1)^2} \norm{y^\perp}_2^4 + \frac{d(d+2)(m-d-1)}{(m-1)^2(m-d-3)} \norm{y^\perp}_2^4 
    \label{eq:approx_error_full}
\end{align}
Now, we look at the following conditional expectation.
\begin{align*}
    \E\left[ (\widehat x_{shr}-x_{LS})^T\Sigma^{-1} (\widehat x_{shr}-x_{LS})\mid S_1A,S_2A \right]
\end{align*}
Since $\widehat x_1$ and $\widehat x_2$ are independent (even when conditioned on $S_1A,S_2A$), we can first take expectation with respect to $\widehat{x}_1$, and using (\ref{eq:approx_error_full}), this becomes
\begin{align*}
    \E\left[ (\widehat x_S-x_{LS})^T\Sigma^{-1} (\widehat x_S-x_{LS}) \mid S_2A \right] - \frac{d^2(d-2)^2}{(m-1)^2} \E\left[ \frac{1}{\widehat x_2^T\Sigma^{-1}\widehat x_2} \right] \\ + \frac{d(d+2)(d-2)^2(m-d-1)}{(m-1)^2(m-d-3)}\E\left[ \frac{1}{\widehat x_2^T\Sigma^{-1}\widehat x_2} \right]
\end{align*}
where $\Sigma=\frac{1}{m}\norm{y^\perp}_2^2 (A^TS_2^TS_2A)^{-1}$ and $\widehat x_S$ (the un-approximated estimator) is calculated using $S_2$ and $\widehat x_2$.
Then we use the normalized error from Lemma \ref{lem:stein} to get 
\begin{align*}
    \E\left[ (\widehat x_S-x_{LS})^T\Sigma^{-1} (\widehat x_S-x_{LS})\mid S_2A \right] &= d - (d-2)^2 \E\left[ \frac{1}{\widehat x_2^T\Sigma^{-1}\widehat x_2} \right]
\end{align*}
Substituting $\Sigma=\frac{1}{m}\norm{y^\perp}_2^2 (A^TS_2^TS_2A)^{-1}$, and using $\E\frac{1}{X}\geq\frac{1}{\E X}$,
\begin{align*}
    &\frac{m}{\norm{y^\perp}_2^2} \E\left[\norm{S_2A(\widehat x_{shr}-x_{LS})}_2^2 \mid S_2A \right] \\
    & = d - \E\left[ \frac{1}{\widehat x_2^T\Sigma^{-1}\widehat x_2} \right]\left( (d-2)^2 + \frac{d^2(d-2)^2}{(m-1)^2} - \frac{d(d+2)(d-2)^2(m-d-1)}{(m-1)^2(m-d-3)} \right)\\
    &\leq d - \frac{(d-2)^2} {d + m\norm{S_2Ax_{LS}}_2^2/ \norm {y^\perp}_2^2} \left( 1 + \frac{d^2}{(m-1)^2} - \frac{d(d+2)(m-d-1)}{(m-1)^2(m-d-3)} \right)
\end{align*}
At this stage we replace the notation $S_2$ back by $S$ which has the identical distribution. Taking expectation with respect to $SA$ on both sides,
\begin{align*}
    &\E\left[\norm{SA(\widehat x_{shr}-x_{LS})}_2^2 \right] \\ &\leq \frac{\norm{y^\perp}_2^2}{m} \E_{SA}\left[ d - \frac{(d-2)^2 \left( 1 + \frac{d^2}{(m-1)^2} - \frac{d(d+2)(m-d-1)}{(m-1)^2(m-d-3)} \right)} {d + m\norm{SAx_{LS}}_2^2/ \norm {y^\perp}_2^2} \right]\\
    &\leq \frac{\norm{y^\perp}_2^2}{m} \left(d - \frac{(d-2)^2 \left( 1 + \frac{d^2}{(m-1)^2} - \frac{d(d+2)(m-d-1)}{(m-1)^2(m-d-3)} \right)} {d + m\E_{SA}[\norm{SAx_{LS}}_2^2]/ \norm {y^\perp}_2^2} \right) & &\text{ from }\E\frac{1}{X}\geq\frac{1}{\E X}\\
    &= \frac{d}{m} \norm{y^\perp}_2^2 \left(1 - \frac{(d-2)^2/d \left( 1 + \frac{d^2}{(m-1)^2} - \frac{d(d+2)(m-d-1)}{(m-1)^2(m-d-3)} \right)} {d + m\norm{Ax_{LS}}_2^2/ \norm {y^\perp}_2^2} \right) & &\text{ since }\E[S^TS]=I \\
    &= \frac{d}{m} \norm{y^\perp}_2^2 \left( 1- \frac{1 - \epsilon'(d,m)} {1 + ({m}/{d}) \norm{Ax_{LS}}_2^2 /\norm{y^\perp}_2^2} \right)
\end{align*}
where $\epsilon'(d,m) = \frac{4(d-1)}{d^2} + \frac{2(d-2)^2}{d(m-1)(m-d-3)}$.

Finally for proving part (c), we use the fact that $S$ is an Oblivious Subspace Embedding \cite{Nelson2016} so that when $m>\mathcal{O}\left( (d+\log\left(1/\delta\right))/\epsilon^2\right)$, with probability greater than $1-\delta$,
\begin{equation*}
(1-\epsilon)\norm{Ax}_2^2 \leq \norm{SAx}_2^2 \leq (1+\epsilon)\norm{Ax}_2^2
\end{equation*}
Then for sufficiently small $\epsilon$ such that $(1-\epsilon)^{-1} \approx(1+\epsilon)$,
\begin{align*}
    \E\left[\norm{A(\widehat{x}_{shr}-x_{LS})}_2^2 \right] &\leq (1+\epsilon)\E\left[ \norm{SA(\widehat{x}_{shr} - x_{LS})}_2^2 \right] \nonumber \\
    &\leq (1+\epsilon) \frac{d}{m} \norm{y^\perp}_2^2 \left( 1- \frac{1 - \epsilon'(d,m)} {1 + ({m}/{d}) \norm{Ax_{LS}}_2^2 /\norm{y^\perp}_2^2} \right)
\end{align*}

\newpage

\end{document}